\documentclass[journal]{IEEEtran}
\usepackage{url}
\usepackage{makeidx}
\usepackage{graphicx}
\usepackage{amsmath}
\usepackage{latexsym}
\usepackage{epsfig}
\usepackage{color}
\usepackage[ruled]{algorithm2e}
\SetKwInput{KwData}{Data}
\usepackage{enumerate}
\usepackage[colorlinks,urlcolor=blue,citecolor=blue,linkcolor=blue]{hyperref}

\usepackage[caption=false,font=normalsize,labelfont=sf,textfont=sf]{subfig}
\usepackage{multirow}

\begin{document}

\title{Contingency-Constrained Unit Commitment With Intervening Time for System Adjustments}
\author{Zhaomiao Guo, Richard Li-Yang Chen~\IEEEmembership{Member,~IEEE}, Neng Fan, and Jean-Paul Watson~\IEEEmembership{Member,~IEEE}
\thanks{
Z Guo is with the Civil and Environmental Engineering Department, University of California, Davis, CA, USA. email: zmguo@ucdavis.edu.}

\thanks{
R. L.-Y. Chen is with the Quantitative Modeling and Analysis Department, Sandia National Laboratories, Livermore, CA 94551, USA. email: rlchen@sandia.gov.}
\thanks{
N. Fan is with the Department of Systems and Industrial Engineering, University of Arizona, Tucson, AZ 85721, USA. e-mail: nfan@email.arizona.edu.}
\thanks{
J.-P. Watson is with the Discrete Math and Optimization Department, Sandia National Laboratories, Albuquerque, NM 87185, USA. e-mail: jwatson@sandia.gov.}

}



\maketitle

\begin{abstract}
The $N$-$1$-$1$ contingency criterion considers the consecutive loss of two components in a power system, with intervening time for system adjustments. In this paper, we consider the problem of optimizing generation unit commitment (UC) while ensuring $N$-$1$-$1$ security. Due to the coupling of time periods associated with consecutive component losses, the resulting problem is a very large-scale mixed-integer linear optimization model. For efficient solution, we introduce a novel branch-and-cut algorithm using a temporally decomposed bilevel separation oracle. The model and algorithm are assessed using multiple IEEE test systems, and a comprehensive analysis is performed to compare system performances across different contingency criteria. Computational results demonstrate the value of considering intervening time for system adjustments in terms of total cost and system robustness.
\end{abstract}


\section*{Nomenclature}
{\noindent Index Sets and Indices}
\begin{itemize}
\item $\mathcal N$: the set of buses, indexed by $i,j$
\item $\mathcal G$: the set of generators, indexed by $g$, $G=|\mathcal G|$
  \begin{itemize}
    \item $\mathcal G_i$: the set of generators located at bus $i$
  \end{itemize}
\item $\mathcal E$: the set of transmission lines, $E=|\mathcal E|$
  \begin{itemize}
    \item $\mathcal E_{.i}$: the set of lines oriented into bus $i$
    \item $\mathcal E_{i.}$: the set of lines oriented out of bus $i$
    \item $(i,j)$: the head bus $i$ and tail bus $j$ of  line $e$
  \end{itemize}
\item $\mathcal{C}$: the set of all $N$-$1$-$1$ contingencies, $C=|\mathcal C|$
\begin{itemize}
\item $c:$ contingency index, $c \in \{1,\ldots,C\}$
\item $\boldsymbol c:$ a binary vector that prescribes a contingency 
\end{itemize}
\item $\mathcal{T}$: the set of time periods, indexed by $t$, $T=|\mathcal T|$
\end{itemize}

{\noindent Parameters}
\begin{itemize}
  \item $B_e$: susceptance  of line $e$
  \item $\overline {\boldsymbol f}$: vector of line capacity
  \item $d_i^{t}$: demand at bus $i$ in period $t$
    \item $\boldsymbol d^{t}$: demand vector in period $t$
  \item $\underline{p}_g$, $\overline{p}_g$: capacity lower/upper bounds for generator $g$
\item $\underline r_g(\boldsymbol x):$ vector of ramp-down rates given unit commitment vector $\boldsymbol x$
\item $\overline r_g(\boldsymbol x):$ vector of ramp-up rates given unit commitment vector $\boldsymbol x$
\item $c^s(\boldsymbol x):$ start-up and shut-down cost given unit commitment vector $\boldsymbol x$
\item $c^p(\boldsymbol p):$ production cost given generation level vector $\boldsymbol p$
  \item $o_e$: allowable line overload factor during secondary contingency periods, $\geq 1$
  \item $\varepsilon$: load shedding threshold, $\varepsilon \in [0,1]$
\end{itemize}

{\noindent Decision Variables}
\begin{itemize}
  \item $\boldsymbol x \in \{0,1\}^{G\times T}$:  unit commitment vector
    \item $\boldsymbol x_g \in \{0,1\}^{T}$:  generator $g$ unit commitment vector
   \item $\boldsymbol p^t, \boldsymbol f^t,\boldsymbol \theta^t$: vectors of generation levels, power flows, and phase angles, respectively, in time period $t$ under the no-contingency (base) scenario
    \item $p_g^t,f_e^t,\theta_i^t$: the generation level of unit $g$, the power flow on transmission line $e$, and the phase angle on bus$i$, respectively, in time period $t$ under the no-contingency (base) scenario
   \item $\boldsymbol p^{ct}, \boldsymbol f^{ct},\boldsymbol \theta^{ct}, \boldsymbol q^{ct}$: vectors of generation levels, power flows, phase angles, and loss-of-load, respectively, in time period $t$ under the no-contingency (base) scenario
  \item $p_g^{ct},f_e^{ct},\theta_i^{ct},q_i^{ct}$: the generation level of unit $g$, the power flow on transmission line $e$, the phase angle on bus $i$, and the loss-of-load at bus $i$, respectively, in time period $t$ under contingency scenario $c$
\end{itemize}

\section{Introduction}
\IEEEPARstart{T}{he} North American Electric Reliability Corporation (NERC) develops and enforces standards to ensure the reliability of power systems in North America. The NERC Transmission Planning Standard \cite{NERC2011} defines system performance requirements under both normal and various contingency conditions. Among contingency conditions, the loss of a single system component ($N$-$1$) and the near simultaneous loss of multiple system components ($N$-$k$) are well studied.

However, a contingency criterion considering non-simultaneous failures of two components
 has not attracted much attention until recently (see \cite{Chatterjee2010,Fan2012}). This contingency criterion, referred to as $N$-$1$-$1$, refers to the consecutive loss of two components with an intervening time period for operator adjustments. In \cite{Chatterjee2010}, the authors performed an $N$-$1$-$1$ contingency analysis of the Midwest ISO's balancing area. More recently \cite{Fan2012} used interdiction methods to study $N$-$1$-$1$ contingency constrained optimal power flows with fixed unit commitment decisions. 
 
Following~\cite{Fan2012}, we assume that an $N$-$1$-$1$ contingency scenario involves the loss of a generating unit or a transmission line, followed by system adjustments (e.g., generator re-dispatch). Following these adjustments, the system experiences a subsequent loss of an additional generator or line. Thus, there are three distinct time periods in a given $N$-$1$-$1$ contingency occurrence. The {\em Base Case} refers to the time periods in which ``the power system is in normal steady-state operation, with all components in service that are expected to be in service". The first loss of a component is referred to as the {\em Primary Contingency/Loss}, while the second loss is referred to as the {\em Secondary Contingency/Loss}. 

In this paper, we study the unit commitment (UC) problem under $N$-$1$-$1$ contingency constraints. Unit commitment involves determination of a minimal-cost ``on-off" schedule of generating units and their respective dispatch levels, subject to physical and operational constraints, in order to satisfy forecasted demand in each time period of the subsequent day. The basic UC problem, discounting contingencies, is well-studied; relevant efforts were reviewed by \cite{Hobbs2001} and more recently by \cite{Padhy2004,Zheng2014}. 

The UC problem with contingency constraints has received increasing attention from academics and practitioners since the 2003 northeast blackout in North America. More specifically, $N$-$1$ and $N$-$k$ contingency-constrained UC and related grid design problems have been studied by \cite{Hedman2010,O'Neill2010} and \cite{Street2011,Wang2012,Chen2014,Chen2015}, respectively.  References \cite{Hedman2010,O'Neill2010} consider $N-1$ contingency constrained UC, and employ line switching to alleviate congestion and yield a more economical dispatch of generation resources. References \cite{Street2011,Wang2012} use robust optimization to find an optimal dispatch schedule under the worst-case $N$-$k$ contingency scenario. Finally, \cite{Chen2014,Chen2015} introduce a new $N$-$k$-$\boldsymbol \varepsilon$ criterion which dictates that at least $(1-\varepsilon^j)$ fraction of the total system demand must be met following the failures of $j$ system components (for $j \in \{,\ldots,k\}$); several decomposition methods were proposed to solve the resulting large-scale optimization model.

In contrast to the above work, we consider $N$-$1$-$1$ contingencies, as defined in \cite{NERC2011}. During the base (no-contingency) case, and during time periods after the primary loss, all thermal limits must be within applicable ratings and loss-of-load is not permitted as a recourse action. During time periods after the secondary loss, controlled load shedding and overloads of transmission lines are allowed for emergency control per NERC standards \cite{NERC2011}. To solve the UC problem with $N$-$1$-$1$ contingency constraints, we introduce a novel decomposition method combining branch-and-cut and a temporal decomposition based on a separation oracle.

The reminder of this paper is organized as follows. In Section \ref{sec:model}, we briefly introduce the baseline UC problem and study the impact of imposing $N$-$1$-$1$ contingency constraints. We then formulate the $N$-$1$-$1$ contingency-constrained UC problem as a large-scale mixed-integer linear program (MILP). In Section \ref{sec:solution}, we describe our novel solution strategy. Numerical experiments on several IEEE test systems are considered in Section \ref{sec:experiments}, where we perform a detailed analysis comparing the impacts of the different contingency criteria on UC solutions. Finally, we conclude in Section \ref{sec:conclusions} with a summary of our results.

\section{Unit Commitment Models}\label{sec:model}
We begin by introducing the baseline unit commitment (BUC) problem. The objective of the BUC problem is to determine a minimal-cost on/off schedule and corresponding dispatch levels for a set of thermal generating units under the no-contingency scenario.  Building on the BUC problem, we then introduce constraints to support $N$-$1$-$1$ contingency compliance. We refer to the extended problem as the $N$-$1$-$1$ CCUC problem, or for conciseness simply $N$-$1$-$1$.

\subsection{The Baseline Unit Commitment Problem}
\label{sec:buc}

The BUC problem described below is based on the deterministic UC formulations introduced in \cite{Carrion2006} and \cite{Wu2010}. We extend these formulations to include a DC approximation of power flow on the transmission network. The BUC problem is formulated as follows:
\begin{subequations}\label{buc}
\begin{align}
\min_{\boldsymbol {x,f,p,\theta}}  \quad &  \boldsymbol c^s(\boldsymbol x) +c^p(\boldsymbol p) \label{mod_uc_obj} \\
\text{s.t.} \quad &\boldsymbol x \in \mathcal X \label{mod_uc_consts}\\
& H \boldsymbol p^t + A \boldsymbol f^t  = \boldsymbol d^t , \ \forall t \label{buc_bal}\\
&	B_{e}(\theta_{i}^{t}-\theta_{j}^{t})- f_{e}^{t}=0,\   \forall e=(i,j), t \label{buc_kirchoff}\\
&|\boldsymbol f^t| \le \overline {\boldsymbol f},\  \forall t\\
	&\underline{p}_g x_g^t\leq p_g^{t}\leq \overline{p}_g x_g^t,\  \forall g, t \\
&\underline{r}_g(\boldsymbol x_g) \le p_g^{t} - p_g^{t-1} \leq \overline{r}_g(\boldsymbol x_g),\  \forall g, t \label{buc_ramp}
\end{align}
\end{subequations}

The objective \eqref{mod_uc_obj} is to minimize the sum of startup and shutdown cost $\boldsymbol c^s(\boldsymbol x)$ and generation cost $c^p(\boldsymbol p)$. With dispatch levels prescribed by $\boldsymbol p$, $c^p (\boldsymbol p)$ is often approximated by a convex quadratic function for thermal generation units. Constraints \eqref{mod_uc_consts} enforce generator minimum uptime/downtime requirements and prescribe startup and shutdown cost as a function of the units committed. The full description of these constraints is provided in Appendix A. Constraints \eqref{buc_bal}-\eqref{buc_ramp} implement economic dispatch under a DC power flow model.  They include (in order): power balance at each bus; power flow on a line; capacity limits for transmission lines; generator dispatch lower and upper bounds; and generator ramping limits across two consecutive time periods. By employing piecewise linearization of the quadratic cost function, the BUC problem \eqref{buc} can be reformulated as a MILP. 

{\bf Remark 1.}  We do not explicitly impose reserve margins in the BUC as $N$-$1$-$1$ compliancy is a stronger reliability requirement; not only does it ensure sufficient generation reserves for all contingency scenarios but additionally considers the placement of these reserves given constraints on transmission availability and capacity (see also  \cite{Hedman2010}).

\subsection{$N$-1-1 Contingency Constrained Operations}
\label{sec:N-1-1}

\subsubsection{Set of all $N$-$1$-$1$ contingencies}
An $N$-$1$-$1$ contingency refers to the loss of two system components in different time periods. We assume losses are possible for any generating unit(s) and/or transmission line(s). Considering all pairs of time periods for primary and secondary losses, and the possible loss of any generating unit and / or transmission line, the set $\mathcal C$ of all $N$-$1$-$1$ contingency scenarios is defined as follows:
\begin{subequations}\label{setD}
\begin{align}
\hskip -0.3cm \mathcal{C} =   \Bigg \{ & \boldsymbol c \in \{0,1\}^{(G+E)\times T} \\
&\sum_{g \in \mathcal G} c_g^t + \sum_{e \in \mathcal E} c_e^t \le 1, \forall t \in \mathcal T \label{cons_at_most_1}\\
&\sum_{t \in \mathcal T} c_g^t \le 1,\quad \forall g \in \mathcal G \label{cons_g1}\\
&\sum_{t \in \mathcal T} c_e^t \le 1,\quad \forall e \in \mathcal E \label{cons_e1}\\
&\sum_{t \in \mathcal T}\sum_{g \in \mathcal G} c_g^t + \sum_{t \in \mathcal T}\sum_{e \in \mathcal E} c_e^t = 2. 
\label{cons_2_exact2} \Bigg \}
\end{align}
\end{subequations}
An $N$-$1$-$1$ contingency scenario $\boldsymbol c \in \mathcal C$ specifies the following: (i) two time periods, denoted by $t_1$ and $t_2$, respectively representing the time periods of primary loss and secondary losses and (ii) one failed component in each of the two periods, denoted by $c_e^{t_1} = 1$ or $c_g^{t_1} = 1$ and $c_e^{t_2} = 1$ or $c_g^{t_2} = 1$. Constraints \eqref{cons_at_most_1} dictate that at most one component can fail in any given time period. Constraints \eqref{cons_g1} and \eqref{cons_e1} specify that each component fails at most once. Constraint \eqref{cons_2_exact2} requires that exactly two distinct components fail. Based on \eqref{setD}, there are ${T\choose 2}=\frac{T(T-1)}{2}$ possible pairs of time periods for primary and secondary losses, $G+E$ possible primary losses, and $G+E-1$ possible secondary losses. Thus, the set $\mathcal C$ has cardinality $|\mathcal C| = C =\frac{T(T-1)}{2}(|G|+|E|)(|G|+|E|-1)$. Clearly, even for moderately-sized power systems and small numbers of time periods, solution of the $N$-$1$-$1$ will pose a considerable computational challenge.


\subsubsection{$N$-$1$-$1$ contingency requirements}
As illustrated in Fig. \ref{fig1}, an $N$-$1$-$1$ contingency scenario is composed of three non-overlapping periods, defined as follows:

\begin{figure}[ht]
\begin{center}
\vskip -0.5cm
\includegraphics[width=0.65\textwidth,angle=0]{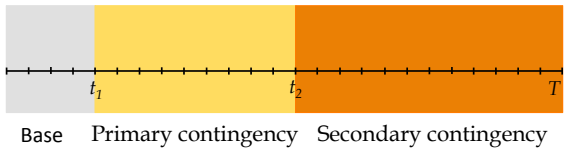}
\end{center}
\vskip -6.7cm
\caption{Three non-overlapping periods of an $N$-$1$-$1$ contingency scenario.}\label{fig1}
\end{figure}

\begin{itemize}
\item \emph{Base} $(t \in \{1,\ldots,t_1-1\}).$ The system operates under normal conditions with no failed components.  Non-anticipativity is enforced during this state, i.e., the operating state (e.g. generation outputs and power flows) are fixed regardless of the specific impeding contingency scenario. 

\item \emph{Primary contingency} $(t \in \{t_1,\ldots,t_2-1\})$. The system operates under a single failed component. At time period $t_1$, the system observes the primary loss and transitions from the nominal operating state (prescribed by $\boldsymbol {p, f, \theta}$) to the contingency $c$ operating state (prescribed by $\boldsymbol p^c, \boldsymbol f^c, \boldsymbol \theta^c$).  If the outaged component is a generator $g$, $p_g^{t_1} = 0$. For all other generators $g'\in \mathcal G \setminus g$, $p_{g'}^{t_1}$ is ramp-constrained by the generator's dispatch level in period $t_1-1$.


\item \emph{Secondary contingency} $(t \in \{t_2,\ldots,T\})$. The system operates with two failed components.  Per NERC reliability standards, controlled load shedding and line overloads are permissible. Therefore, the line overload factor $o_e$ and allowable load shedding $\varepsilon$ (as a fraction of total demand) can be utilized to alleviate operational infeasibilities. 

\end{itemize}

For conciseness, we introduce in the following "in the contingency'' $\boldsymbol w$ indicators and both allowable loss-of-load $\boldsymbol h$ and line overload $\boldsymbol o$ quantities:
\begin{align}\sum_{t = 1}^t c_e^{ct} = w_e^{ct}, \ \forall e, t \text{ and }\sum_{t = t}^t c_g^{ct} = w_g^{ct}, \ \forall g, t
\end{align}
%
\begin{align}
h_i^{ct} =
\begin{cases}
    0, \quad &\forall i, t = 1,\cdots,t_2^c -1\\
    d_i^t, \quad &\forall i, t = t_2^c,\cdots,T
\end{cases}\label{uc-cons2}
\end{align}

\begin{align}
o_e^{ct} =
\begin{cases}\label{uc-cons3}
    0, \quad &\forall e, t = 1,\cdots,t_2^c -1 \\
    o_e, \quad &\forall  e, t = t_2^c,\cdots,T
\end{cases}
\end{align}

Constraints \eqref{uc-cons2} and \eqref{uc-cons3} require that during secondary contingency periods, the allowable loss-of-load and overload factor are equal to $d_i^t$ and $o_e$, respectively.  These values are zero in all other periods.

Then, $\forall \boldsymbol c \in \mathcal C$, the contingency dispatch operation in periods $t \in \{t_1^c,\ldots, T\}$ is constrained as follows:
\begin{subequations}\label{ccons}
\begin{align}
& H \boldsymbol p^{ct} + A \boldsymbol f^{ct} + \boldsymbol q^{ct} = \boldsymbol d^t, \ \forall t \label{cc_bal}\\
& B_{e}(\theta_{i}^{ct}-\theta_{j}^{ct})(1-w_e^{ct}) - f_{e}^{ct}=0,\   \forall e=(i,j), t \label{cc_vpa}\\
&|f_e^{ct}| \le \overline f_e(1-w_e^{ct})(1+o_e^{ct}),\  \forall e,t \label{cc_line_cap}\\
&\underline{p}_g x_g^t(1-w_g^{ct})\leq p_g^{ct}\leq \overline{p}_g x_g^t(1-w_g^{ct}),\  \forall g, t \label{cc_gen_cap}\\
&\underline{r}_g(\boldsymbol x_g) \le p_g^{ct} - p_g^{c,t-1} \leq \overline{r}_g(\boldsymbol x_g),\  \forall g, t \label{cc_gen_ramp}\\
& \boldsymbol q^{ct} \le \boldsymbol h^{ct}, \ \forall t \label{cc_loss_ub}\\
& \boldsymbol 1^\top \boldsymbol q^{ct} \le \varepsilon \boldsymbol 1^\top \boldsymbol h^{ct}, \ \forall t \label{cc_tot_loss}\\
&p_g^{ct_1^c-1} - {p_g^{t_1^c - 1}}=0, \ \forall g \in \{g | c_g^{t_1^c} = 0\} \label{cc_link}
\end{align}
\end{subequations}

Contingency constraints include (in order of appearance): power balance at each bus, with unsatisfied demand $\boldsymbol q^{ct}$ \eqref{cc_bal}; power flows on each line \eqref{cc_vpa}; line capacity bounds \eqref{cc_line_cap} with overload factor $o_e$ during secondary contingency periods; generation dispatch bounds \eqref{cc_gen_cap}; generator ramping limits \eqref{cc_gen_ramp}; upper bound on loss-of-load at each bus \eqref{cc_loss_ub}; threshold for total loss-of-load in periods after secondary loss \eqref{cc_tot_loss}; and the non-anticipativity constraint for generators not in the contingency in period $t_1^c$ \eqref{cc_link}.


\subsection{Full Formulation}

The optimization objective in the $N$-$1$-$1$ model is to find a minimum-cost UC and no-contingency scenario economic dispatch such that a feasible recourse power flow exists under each $N$-$1$-$1$ contingency scenario.  The full formulation is obtained by combining the BUC model \eqref{buc} with the full set of contingency constraints \eqref{ccons}, one for each contingency scenario $\boldsymbol c \in \mathcal C$. The full formulation of the $N$-$1$-$1$ CCUC model is then given as follows.
\begin{subequations}\label{full_mod}
\begin{align}
\min_{\substack{\boldsymbol{x,f,p,\theta} \\ \boldsymbol f^c, \boldsymbol p^c, \boldsymbol q^c, \boldsymbol \theta^c}}\ &c^s(\boldsymbol x) + c^p(\boldsymbol p) \label{full_mod_obj} \\
\text{s.t.}\quad & \text{Constraints }\text{\eqref{mod_uc_consts} -- \eqref{buc_ramp}} \label{full_mod_uc_cons}\\
&\text{Constraints } \text{\eqref{cc_bal} -- \eqref{cc_link}},\ \forall \boldsymbol c \in\mathcal{C} \label{full_mod_ccons}
\end{align}
\end{subequations}

Model \eqref{full_mod} is an extremely large-scale MILP due to the full set of DCOPF constraints \eqref{full_mod_ccons}, one for each contingency scenario.  The objective \eqref{full_mod_obj} includes only the unit commitment cost and the no-contingency scenario generation cost.  However, extension to consider worst-case cost is straightforward. Following established models (\cite{Hedman2010}, \cite{Chen2014}, \cite{Chen2015}), we ignore costs during a contingency state, as the primary concern of system operators during a contingency is to ensure operational feasibility and system stability.  


{\bf Remark 2.} In the $N$-$1$-$1$ model, there are $\frac{T(T-1)}{2}(G+E)(G+E-1)$ sets of constraints \eqref{full_mod_ccons}, which collectively ensure that a feasible recourse power flow exists in each contingency scenario. When defining the full set of contingency scenarios \eqref{setD}, we assumed that when a primary contingency component fails, its failure persist for the remainder of the planning horizon. We can relax this assumption through introduction of a new integer parameter $\tau \ge 1$ that prescribes the number of time periods until the primary contingency component is returned to service. Under this assumption we do not need to consider all time period pairs for the primary and secondary contingency. Rather, it then suffices to consider all pairs of time periods whose difference is $\tau$. Then, there are $(T-1)+(T-2)+\cdots+(T-\tau)=\frac{(2T-\tau-1)\tau}{2}$ pairs of primary and secondary failure periods $t_1$ and $t_2$. 

%

\section{Solution Approaches}
\label{sec:solution}

The full MILP model \eqref{full_mod} can be solved using an iterative algorithm like Benders decomposition; direct solution via the extensive form is not practical. In applying BD, we first decompose the full problem into a master problem (MP), defined by \eqref{full_mod_obj}-\eqref{full_mod_uc_cons},  and a set of subproblems (SP), defined by \eqref{full_mod_ccons}, one for each contingency $\boldsymbol c \in \mathcal C$. The MP prescribes the unit commitment vector $\boldsymbol x$ and the no-contingency economic dispatch $(\boldsymbol {f, p, \theta})$. The subproblems SP$(\boldsymbol {x, p, c})$ are based on Constraint  \eqref{ccons}, with the following augmentations:
\begin{itemize}
\item Add variables $\boldsymbol s$ to indicate the amount of load shedding above the allowable threshold $\varepsilon$. This ensures relatively complete recourse.
\item Add the trivial objective function $\min \ \boldsymbol 1^\top \boldsymbol s$. 
\end{itemize}
 The objective of SP is to minimize the amount of demand shed above the allowable threshold. If the objective value of SP$(\boldsymbol {x, p, c})$, given by $z$, is zero then the current solution $\boldsymbol {x, p}$ can survive contingency scenario $\boldsymbol c$. Otherwise, a Benders feasibility cut can be added to the MP to eliminate the infeasible solution $(\boldsymbol {x,p}).$

In state-of-the-art algorithms for contingency-constrained UC problems, BD is typically implemented as a cutting plane algorithm (CPA) (\cite{Chen2014},\cite{Chen2015},\cite{Wu2010},\cite{Street2011},\cite{Yao2007}, and \cite{Yuan2014}). 
Although easy to implement, CPA has limited computational tractability because at each iteration of the algorithm an integer program MP must be solved to select a candidate UC decision. We next outline a branch-and-cut algorithm (BCA) that avoids this drawback by only exploring the branch-and-bound tree corresponding to the UC $\boldsymbol x$ once. BCA is a branch-and-bound algorithm in which cutting planes are generated within the branch-and-bound tree. 
%

\subsection{Branch-and-Cut Algorithm}

Recent advances in optimization solver technology, both in commercial (IBM CPLEX \cite{cplex} and Gurobi \cite{gurobi}) and non-commercial (SCIP \cite{scip} and DIPS \cite{dip}) packages, have significantly simplified the implementation of branch-and-cut algorithms by enabling the addition of violated inequalities directly within the branch-and-bound tree, thus avoiding the need to repeatedly explore the branch-and-bound tree defined by the binary UC variables $\boldsymbol x$. In IBM CPLEX, branch-and-cut algorithms can be implemented using IloCplex.Callback functions (e.g., IloCplex.LazyConstraintCallback).  In recent years, work on a number of combinatorial problems (e.g., survivable network design  \cite{Fortz2009} and the minimum tool booth \cite{Bai2009}) have shown that a significant reduction, often an order of magnitude or better, in computational time can be achieve using a BCA compared to a CPA, with increasing reductions in runtime for harder and larger combinatorial problems.  This motivates our development of a branch-and-cut algorithm for $N$-$1$-$1$. 

Let the linear programming relaxation of MP, the \emph{relaxed master problem} (RMP), be given as follows:
\begin{subequations}\label{rmp}
\begin{align}
\min_{\boldsymbol {x, f, p, \theta}} \quad  &c^s(\boldsymbol x) + c^p(\boldsymbol p) \\
\textmd{s.t.} \quad & \boldsymbol x \in \mathcal X\\
&\boldsymbol 0 \le \boldsymbol x \le \boldsymbol 1.
\end{align}
\end{subequations}

Let $L$ be the list of nodes of the B\&B tree to explore. Initially, $L$ contains only the root node $o$ with no branching constraints. Before initializing the BCA, a number of valid inequalities are added to the initial RMP to strengthen the root node LP relaxation.  We refer to these valid inequalities as the cut pool \emph{P}.  We note that in BCA implementations it is important to strengthen the root node adequately to avoid unnecessary exploration of the B\&B tree.  Starting with an empty \emph{P} results in very slow pruning because many infeasible nodes will  only be pruned late in the search.  However, adding too many valid inequalities to \emph{P} slows down the LP relaxation solves at each node.  Thus, there is a significant trade-off between strengthening the root node LP relaxation, and thus avoiding unnecessary exploration of the B\&B tree, and overburdening the RMP, and thus increasing LP solve times. This trade-off is more art than science and may be very specific to the application domain. In our BCA implementation, we initialize $P$ with the following valid inequalities.

Let $\lambda_1$ and $\lambda_2$ be the capacity of the smallest and the second smallest generators, respectively. Then the following are valid inequalities:
\begin{subequations}
\begin{align}
&\sum_{g \in G} x_g^t \ge 3, \quad \forall t=2,\cdots,T  \label{cons_min_gen_on}\\
&\sum_{g \in G} \overline p_g x_g^1 - \lambda_1 \ge \boldsymbol 1^\top \boldsymbol d^1  \label{cons_min_gen_cap1}\\
&\sum_{g \in G} \overline p_g x_g^t - \lambda_2 \ge \boldsymbol 1^\top \boldsymbol d^t,  \quad \forall t=2,\cdots,T  \label{cons_min_gen_cap2}   
\end{align}
\end{subequations}
Constraints \eqref{cons_min_gen_on} require that at least three generators be on at any given time period after period one, since  two generators may fail in consecutive periods. Constraint \eqref{cons_min_gen_cap1} dictates that in period 1 the maximum capacity across all committed generators minus the capacity of smallest generator in the system must be greater than the total demand of period 1. Constraints \eqref{cons_min_gen_cap2} require the total capacity across of committed generators minus the aggregate capacity of the two smallest generators must be larger than the demand, for each period after the first period. We have considered only a simple set of valid inequalities for $P$ and significant ``tuning'' is required to determine the "optimal" valid inequalities to initialize $P$. 
%
%
%
%

Solving the RMP for a node $o$ means solving the RMP with associated branching constrained defined by $o$. Node $o$ prescribes the subset of discrete variables that are fixed at that particular node of the B\&B tree. In a depth-first variant of the BCA, branching is performed until an integer solution is identified.  Only when an integer solution is identified are violated inequalities screened for each contingency scenario.  At the opposite extreme, violated inequalities for each contingency scenario may be screened after each node solve, resulting in the breath-first variant of the BCA.  The former, depth first variant, may result in a large node list $L$ and the later may result in the addition of a large number of violated inequalities. A compromise between these two extremes may be achieved by defining a branching depth parameter to control the trade-off between branching and cut generation.

\subsection{Temporally Decomposed Bilevel Separation Oracle}

Typically, the existence of a feasible DCOPF for each contingency scenario must be verified explicitly. Such explicit enumeration, however, cannot scale as the number of contingencies $C$ is extremely large, even for moderately-size power systems and/or planning horizons. As we demonstrate in our computational experiments, securing the system against a small number of contingency scenarios is empirically sufficient to ensure feasibility against the full set of contingency scenarios $\mathcal C$. So, instead of explicitly checking feasibility across all contingency scenarios, we describe a bilevel separation problem to screen  a small number of worst-case contingency scenarios that in aggregate ``covers'' $\mathcal C$.

As noted in \cite{Chen2014} and \cite{Chen2015}, significant computational challenges exist in solving bilevel programs for power system vulnerability analysis. In $N$-$1$-$1$ UC, these challenges are further compounded by the fact that the vulnerability analysis problem spans the $T$ planning periods, resulting in a very large-scale bilevel program with $T\times G\times E$ upper-level binary decision variables. To overcome this computational challenge, we perform a temporal decomposition in which component failures are restricted to a preselected time period pair.  We then iterate over all possible time periods pairs, solving a simplified and smaller bilevel separation oracle for each time period pair.

Given a unit commitment schedule $\boldsymbol {x}$, a no-contingency scenario power flow $ {\boldsymbol {\tilde p}}$, and a pair of periods $t_1$ and $t_2$ denoting the times of the primary and the secondary contingency, respectively, the \emph{Power System Inhibition Problem} (PSIP) can be used  to determine the worst-case loss-of-load under any two non-simultaneous component failures. In this bilevel program, the upper-level decision vector $\boldsymbol c$ identifies the primary and secondary component losses, and the lower-level decision vectors $(\boldsymbol {f,p,q,s, \theta})$ correspond to recourse DCOPF under the contingency scenario prescribed by $\boldsymbol c$.

The set of valid contingency scenarios given a pair of time periods $(t_1, t_2)$ is given as follows:
\begin{subequations}\label{setD2}
\begin{align}
\mathcal{C}(t_1,t_2) =   \Bigg \{ & \ \  \boldsymbol c \in \{0,1\}^{(G+E)\times T}: \\
&\sum_{g \in G} c_g^t + \sum_{e \in \mathcal E} c_e^t = 
\begin{cases}\label{D2_two_failures}
    1,\ \forall t \in \{t_1, t_2\} & \\
    0,\ \text{otherwise},&\\
\end{cases} \\
&c_g^{t_1} + c_g^{t_2} \le 1,\quad \forall g \in \mathcal G, \label{D2_single_failure_g}\\
&c_e^{t_1} + c_e^{t_2} \le 1,\quad \forall e \in \mathcal E \label{D2_single_failure_e} \ \
\Bigg \}
\end{align}
\end{subequations}

Constraint \eqref{D2_two_failures} dictates that a single component failure occurs in each of time periods $t_1$ and $t_2$. Constraints \eqref{D2_single_failure_g} and \eqref{D2_single_failure_e} require that each component can fail at most once. 

In preceding formulations, $\boldsymbol c$ was an input parameter. However, in the PSIP $\boldsymbol c$ is decision vector prescribing the two failed components.  For notational convenience, we again define binary variables $\boldsymbol w$ to indicate whether or not a component failed in the current period or a prior period. Then PSIP$(\boldsymbol x, {\boldsymbol {\tilde p}}, t_1, t_2)$ is given as follows:
\begin{subequations}\label{bl_psip}
\begin{align}
\max_{  \substack{  \boldsymbol w\\ \boldsymbol c \in \mathcal C(t_1, t_2)   }  }\ &\min_{\boldsymbol {f,p,q,s,\theta}} \quad  \boldsymbol 1^\top \boldsymbol s \label{bl_psip_obj}\\
\textmd{s.t.}\quad & H \boldsymbol p^t + A \boldsymbol f^t + \boldsymbol q^t= \boldsymbol d^t, \ \forall t\\
&B_{e}(\theta_{i}^{t} - \theta_{j}^{t}) (1-w_e^t) - f_{e}^{t} = 0, \   \forall e,t\\
&|f_e^t| \le \overline f_e^t ( {1-w_e^t})(1-o_e^t), \   \forall e,t\\
&\underline{p}_g x_g^t (1-w_g^t) \le p_g^t \le \overline{p}_g^t x_g^t (1-w_g^t), \ \forall g, t\\
&\underline{r}_g(\boldsymbol x) \le p_g^{t} - p_g^{t-1} \leq \overline{r}_g(\boldsymbol x), \ \forall g, t\\
& \boldsymbol 1^\top \boldsymbol q^t - s^t \le \varepsilon \boldsymbol 1^\top \boldsymbol d^t, \ \forall t \\
&p_g^{{t_1}-1} - {   \tilde p_g^{t_1 - 1}}(1-w_g^{t_1})=0, \ \forall g \label{bl_psip_non_anticipativity}\\
&\boldsymbol {p,q,s} \ge \boldsymbol 0
\end{align}
\end{subequations}

The optimization objective \eqref{bl_psip_obj} is to maximize the minimum amount of load shed above the allowable threshold given by $\boldsymbol s$. Constraint \eqref{bl_psip_non_anticipativity} enforces non-anticipativity and ramp constraints any generator not in the contingency at time $p_g^{t_1}$ base on its no-contingency scenario output $\tilde p_g^{t_1-1}$.

Model \eqref{bl_psip} is a bilevel program that cannot be solved directly by commercial solvers, but it can be reformulated as a mixed integer linear program (MILP) by dualizing the lower-level linear program, combing the result with the upper maximization problem, and linearizing the resulting bilinear terms in the objective. In the following algorithmic description, the solution of PSIP refers to the solution of the linearized MILP counterpart of \eqref{bl_psip}.

Instead of solving $C$ linear programs \eqref{ccons} at each iteration, we can solve $\frac{1}{2}(T)(T-1)$ mixed-integer linear programs \eqref{bl_psip} for every time period pairs $(t_1,t_2)$. We refer to this variant algorithm as ``Hybrid Branch-and-Cut'' and present the full algorithm as follows:

\begin{algorithm}\label{hbca}
\caption{Hybrid Branch-and-Cut Algorithm (HBC)}
\textbf{Require:} A starting cut pool \emph{P}\;
 \textbf{Initialization:} $z^* = +\infty, L= \{o\}$ \tcc*[r]{ {\footnotesize $o$ has no branching constraints}}
 \While{L \normalfont{is nonempty}}{
  Select a node $ \hat o \in L$\;
  {$L \gets L \setminus \{\hat o\}$} {\footnotesize { \tcc*[r]{delete node}}}
  Solve \text{RMP} for node $\hat o$. Let $ {(\boldsymbol {x,f,p,\theta})}$ be the optimal solution and $z$ be the optimal cost\;
    \If{$ z \ge z^*$}{
   \eIf{$ {\boldsymbol x}$ \normalfont{is fractional}}{
   Branch, resulting in nodes $o'$ and $o''$, $L \gets L \cup \{o',o''\}$ {\footnotesize { \tcc*[r]{add nodes}}}
   }{
   \ForEach{$(t_1,t_2)$ \normalfont{pair}}{Solve \textsc{psip}$( {\boldsymbol x}, {\boldsymbol p},  t_1, t_2)$. Let $\omega$ be the optimal objective value\;
   
   \eIf{$ \omega>0$}{
   add  $f-cut$ to $P$\;
   $L \gets L \cup \{\hat o\}$ {\footnotesize { \tcc*[r]{put back node}}}}{
\If{$  z < z^*$}{Define new upper bound $z^* \gets z$ and 
   update incumbent $\boldsymbol x^* \gets {\boldsymbol x}$ }
  }}}}}
  \Return $\boldsymbol x^*$
\end{algorithm}

\section{Numerical Experiments and Analysis}
\label{sec:experiments}

We tested our models and the proposed HBC algorithms on the IEEE 6-bus, 9-bus, 14-bus,  24-bus, 30-bus and 39-bus systems \footnote{All test cases can be obtained from http://pserc.cornell.edu/matpower}  \cite{Matpower} on a laptop running Mac OSX, with 2.3 GHz Intel Core i7 processor and 8GB of memory. The models and algorithms are implemented in C++ using IBM's Concert Technology Library 2.9 and CPLEX 12.6 MILP solver. For each instance, we considered 12 time periods and load-shedding threshold $\varepsilon$ and line overload factor $o_e$ of 0.15. The peak load factor for each system was estimated using the Electric Reliability Council of Texas (ERCOT) demand data \cite{ERCOT}. 

\begin{figure}[t]
\begin{center}
\includegraphics[width=2in]{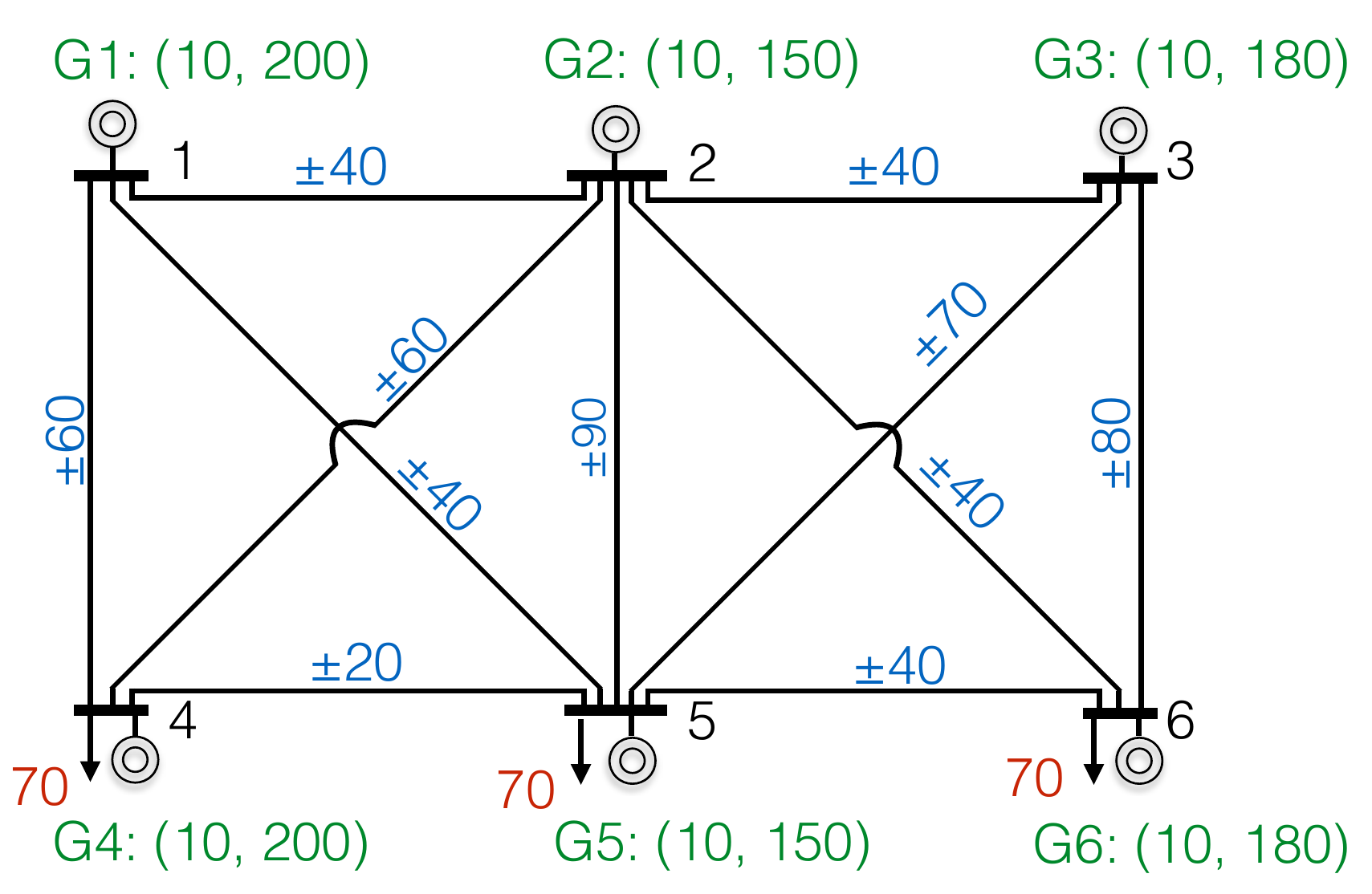}
\end{center}
\vspace{-0.4cm}
\caption{Single-line diagram of the modified IEEE 6-bus test system.}\label{fig: 6_bus_data}
\end{figure}

We begin with an analysis of optimal commitment and dispatch decisions for the 6-bus system across the different security criteria. The single-line diagram of the 6-bus system is shown in Figure \ref{fig: 6_bus_data}, in which green text denotes generator data, with minimum and maximum generation limits in parentheses; blue text denotes maximum power flow in the nominal operating state; and red text denotes loads. We added three fast-ramping (and more costly) units  $G4$ - $G6$ in close proximity to load buses, to ensure system feasibility across the four security criteria:  $N$-$0$ (with the ``1.5 rule'' for reserve margins\footnote{The ``1.5 rule'' is defined as follows: A power system will carry reserves equal to the dispatch level of the largest generator plus one half the dispatch level of the second largest generator. This rule is adopted by several ISOs (e.g., ISO-NE) in the United States}), $N$-$1$, $N$-$1$-$1$, and $N$-$2$.

\begin{figure}[!t]
\begin{center}
\includegraphics[width=2.5in]{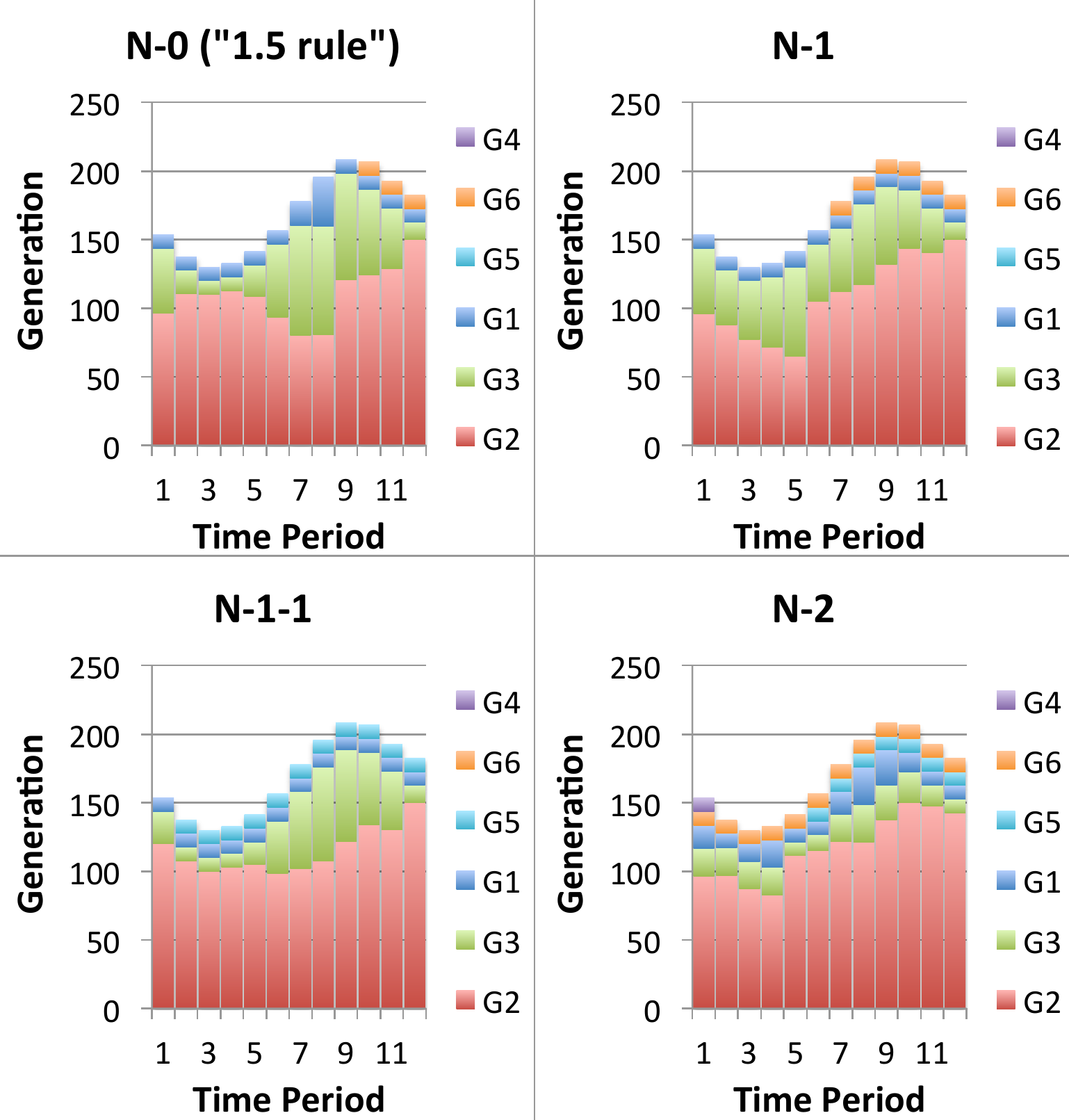}
\end{center}
\vspace{-0.4cm}
\caption{Optimal dispatch levels for different contingency criteria (6-bus)}\label{fig: 6_bus_gen}
\end{figure}

The optimal commitment and dispatch decisions for the 6-bus system varies across different contingency criteria, as shown in see Figure \ref{fig: 6_bus_gen}. For clarity of exposition, we summarize the optimal commitment decisions for  two representative time periods (off-peak and on-peak) in Table \ref{tbl: uc_diff_types}, in order to illustrate the differences between UC decisions under different security criteria.  Under the ``1.5 rule'' ($N$-$0$), only slower-ramping and lower-cost units $G1$ - $G3$ are committed in both time periods. Under the $N$-$1$ criterion, unit $G6$ is committed in addition to $G1$ - $G3$ during the peak period. Under the $N$-$1$-$1$ criterion, we additionally commit $G5$ during the off-peak period and switch from $G6$ to $G5$ during the peak period. The switch from $G6$ to $G5$ reflects a trade-off between production cost and start-up cost. Unit $G6$ has a lower start-up cost but higher production cost than $G5$. Consequently, when a generator may be needed for a longer time (e.g., under $N$-$1$-$1$), $G5$ is preferable. Lastly, when the security criterion includes consideration for two (near) simultaneous failures ($N$-$2$) with no intervening time for adjustments, the additional fast-ramping generator $G6$ is committed. These results are consistent with our intuition that fast-ramping units $G4$ - $G6$ are required only if $G1$ - $G3$ are not able to meet the security requirements. Additionally, comparing different security criteria, $N$-$2$ requires more units to be committed than $N$-$1$-$1$ due to the lack of adjustment time between failures, and more units are committed for $N$-$1$-$1$ than $N$-$1$ in order to cover the loss of a second unit in subsequent periods after the primary loss.

\begin{table}[ht]
\renewcommand{\arraystretch}{1.3}
\caption{UC decisions under different contingent criteria for off-peak and on-peak periods (modified 6-bus system)}
\label{tbl: uc_diff_types}
\centering
\resizebox{\columnwidth}{!}{
\begin{tabular}{c|c c c c c c| c c c c c c}
\hline
\multicolumn{1}{c|}{} & \multicolumn{6}{c|}{Off-peak ($t = 3$)}&\multicolumn{6}{c}{On-Peak ($t = 9$)} \\
 & G1 & G2 &G3 &G4 &G5& G6& G1 & G2 &G3 &G4 &G5& G6\\
\hline
$N-0 $ (``1.5 rule'') &$\surd$  &$\surd$&$\surd$&&& &$\surd$&$\surd$&$\surd$&&\\
$N-1$ & $\surd$&$\surd$&$\surd$&&& & $\surd$&$\surd$&$\surd$&&&$\surd$\\
$N-1-1$ &$\surd$ &$\surd$&$\surd$&&$\surd$& &$\surd$ &$\surd$&$\surd$&&$\surd$&\\
$N-2$ &$\surd$ &$\surd$&$\surd$&&$\surd$&$\surd$ &$\surd$ &$\surd$&$\surd$&&$\surd$&$\surd$\\
\hline
\end{tabular}
}
\end{table}

\begin{figure}[t]
\begin{center}
\includegraphics[width=3in]{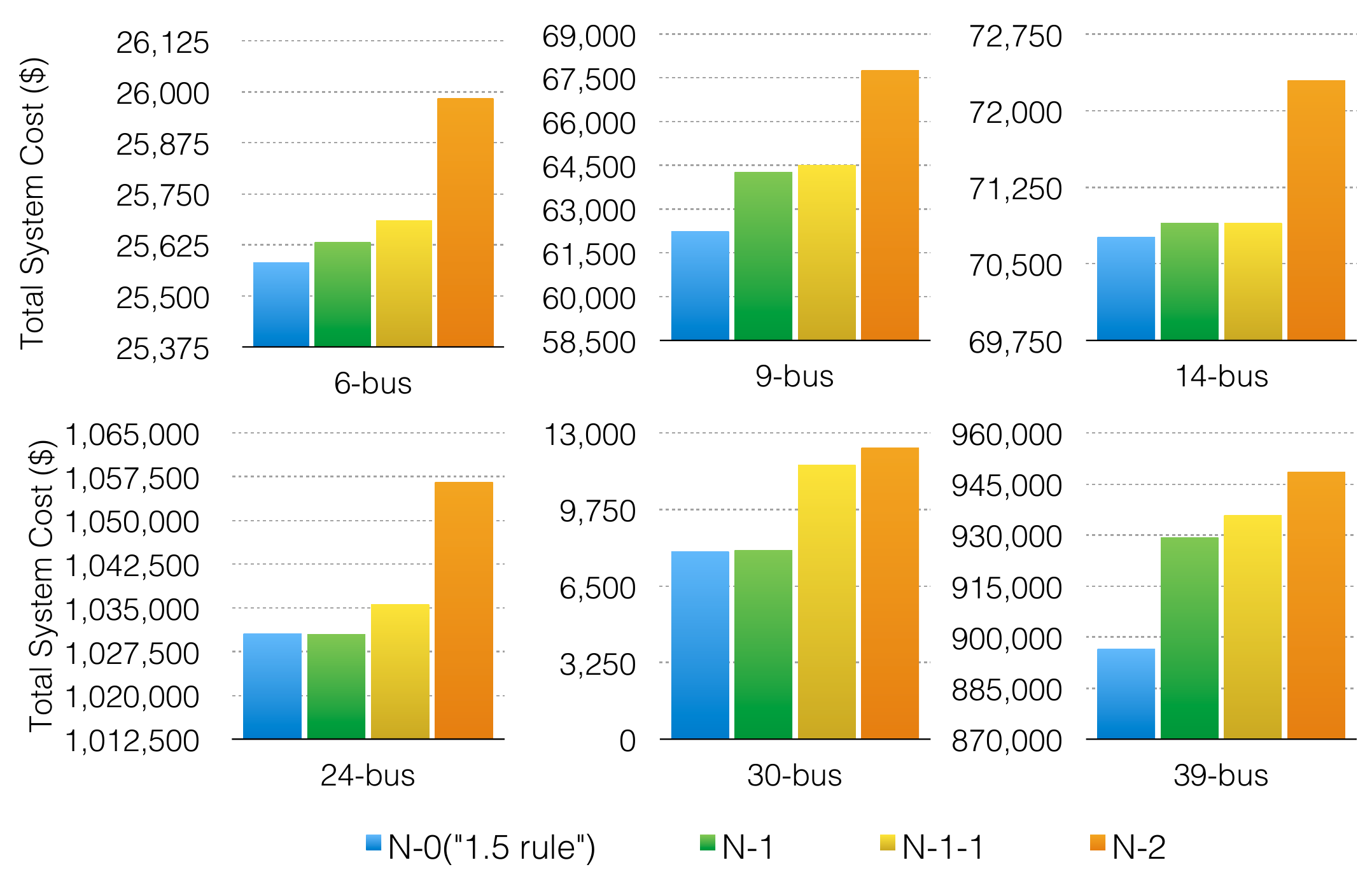}
\end{center}
\vspace{-0.4cm}
\caption{Total costs across different contingency criteria and test systems}\label{fig: obj_value_6_9}
\end{figure}

Next, we compare total costs under different security criteria and test systems, illustrated in Figure \ref{fig: obj_value_6_9}. As expected, we observe that total cost increases monotonically from $N$-$0$ (``1.5 rule''), $N$-$1$, and $N$-$1$-$1$ to $N$-$2$. In all cases, the cost difference under the $N$-$2$ and $N$-$1$-$1$ criteria demonstrate the value of intervening time for system adjustments during a multiple-failure contingency scenario. However, the magnitude of these cost differences are system specific; some test systems (e.g., the 6-bus system) exhibit significant cost differences under the $N$-$1$-$1$ and $N$-$2$ criteria but show little cost difference between $N$-$1$ and $N$-$1$-$1$.  In contrast, the 30-bus system shows the opposite trend: little difference under the $N$-$1$-$1$ and $N$-$2$ criterion but significant  differences between $N$-$1$ and $N$-$1$-$1$.

\begin{figure}[!t]
\begin{center}
\includegraphics[width=2in]{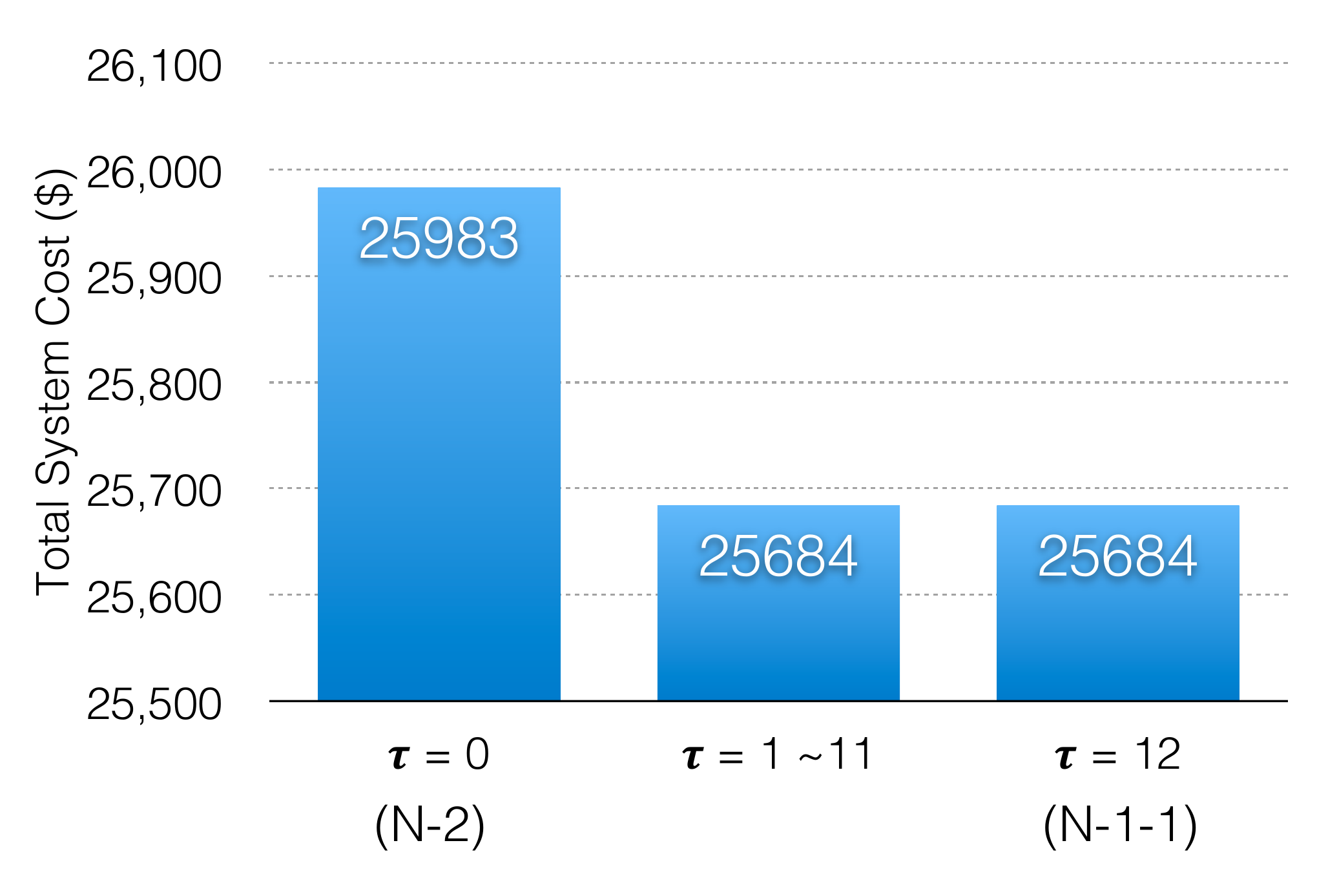}
\end{center}
\vskip -0.5cm
\caption{Total cost as a function of $\tau$}
\label{fig: tau_obj_val2}
\end{figure}

In the context of {\bf Remark 2}, we assess the impact of varying $\tau$, the maximum time difference between the first and the second component failures, for all six test systems. Notice that $\tau=0$ under $N$-$1$-$1$ corresponds precisely to $N$-$2$. Results for the 6-bus system are summarized in Figure \ref{fig: tau_obj_val2}. Qualitatively similar results hold for the other test systems. The results indicate that total cost is invariant to changes in $\tau$ greater than zero. We observe that the computational tractability of the $N$-$1$-$1$ problem is strongly dependent on the parameter $\tau$. As $\tau$ increases, the number of time period pairs ($t_1, t_2$) grows exponentially. Our results suggest that solving $N$-$1$-$1$ instances with $\tau=1$, considering only intervening time of a single period, may be a good proxy for the full $N$-$1$-$1$ problem. However, rigorously mathematical proof and/or comprehensive numerical experiments on a variety of test systems are required to validate this numerical observation. 

%

Finally, we discuss the efficacy of our HBC algorithm and the implicit search for the worst-case contingency scenario by solving PSIPs, rather than explicitly screening the full set of contingency scenarios for each candidate UC solution. We denote the size of the dynamic contingency list as $M$. The size of the final dynamic contingency list $M$ and the cardinality of the contingency scenario set $C$ are reported in Table \ref{tbl: num_fail}. For each test system, the number of contingency scenarios identified is only a very small fraction of the total possible contingency scenarios, which implies that we only need to generate feasibility cuts using a small set of ``sufficient'' contingency scenarios to ensure $N$-$1$-$1$ compliancy. This observation is consistent with the analogous observation for general $N$-$k$ cases, as reported in \cite{Chen2014} and \cite{Chen2015}.

\begin{table}[!t]
\renewcommand{\arraystretch}{1.3}
\caption{Percentage of contingency scenarios identified for the $N$-$1$-$1$ contingency criterion}
\label{tbl: num_fail}
\centering
\begin{tabular}{c c c| c l l}
\hline
Test systems & $G$ & $E$ & $M$ & $C$ & Ratio\\
\hline
6-bus & 6 & 11 & 5 & 330 & 1.52\% \\
9-bus & 6 & 9 & 8 & 330 & 2.42\% \\
14-bus & 10 & 20 & 3 & 990 & 0.30\% \\
24-bus & 32 & 38 & 3 & 10912 & 0.03\% \\
30-bus & 8 & 41 & 10 & 616 & 1.62\% \\
39-bus & 14 & 46 & 16 & 2002 & 0.80\% \\
\hline
\end{tabular}
\end{table}

\section{Conclusions}
\label{sec:conclusions}

We focus on a UC model that considers non-simultaneous component failures with intervening time for system adjustments. The naive formulation of this model is an extremely large-scale MILP, due to the coupling of consecutive time period pairs. To overcome these computational challenges, we introduce an efficient BCA algorithm using a temporally decomposed separation oracle. The model and algorithm are tested on multiple IEEE test systems, through which we observe that (1) the often significant changes in commitment and dispatch decisions for different security criteria; (2) the cost benefit of system adjustment time; (3) the general $N$-$1$-$1$ problem can be well approximated by only considering consecutive failures (i.e., $\tau = 1$); and (4) the number of contingency scenarios identified by our approach represents a very small fraction of the total number.

Allowing for non-simultaneous failures in unit commitment significantly increases the computation burden, relative to other contingency criteria. Although our proposed CPA framework is tractable for small and moderately size test systems, significant computational challenges remain for larger industrial systems and further algorithmic research in required to achieve scalability. A stronger PSIP formulation or an efficient heuristic will be required to cope with the computational challenge posed by large-scale bilevel programs. Finally, a rigorous mathematic proof and/or a comprehensive computational study is needed to support the idea of considering only one single period between successive failures (i.e., $\tau=1$) as a proxy for the fully general $N-1-1$ problem. Aside from computational challenges, from the theoretical point of view the relationship (e.g., the relative magnitude of cost difference) among different contingency criteria is also an interesting area for further research. 

\section{Acknowledgements}

Sandia National Laboratories is a multi-program laboratory managed and operated by Sandia Corporation, a wholly owned subsidiary of Lockheed Martin Corporation, for the U.S. Department of Energy's National Nuclear Security Administration under contract DE-AC04-94AL85000.

\appendices
\section{Full Description of Constraints (\ref{mod_uc_consts})}

The explicit description of Constraints (\ref{mod_uc_consts})\footnote{$T_g^{u0}$/$T_g^{d0}$ denote initial minimal online/offline times; $T_g^{u}$/$T_g^{d}$ denote nominal minimum online/offline times; $C_g^u$/$C_g^d$ denote startup/shutdown costs.} is as follows, which include (in order): initial online and offline requirements for generators; minimum uptime in nominal time periods; minimum uptime for the last $T_g^u$ periods; minimum downtime in nominal time periods; minimum downtime for the last $T_g^u$ periods; startup costs; shutdown costs; non-negativity for startup/shutdown costs; and binary constraints for the on/off status of generators.

\begin{align}
\begin{cases}\label{uc-cons1}
    \sum_{t=1}^{T_g^{u0}} (1-x_g^t) = 0,\ \forall g\\
    \sum_{t=1}^{T_g^{d0}} x_g^t = 0,\ \forall g\\
	\sum_{t'=t}^{t+T_g^u -1} x_g^{t'}\geq T_g^u (x_g^t- x_g^{t-1}),\\
    \hskip 0.7cm  \forall g, t = T_g^{u0}+1, \cdots, T-T_g^u+1\\
	\sum_{t'=t}^T (x_g^{t'}-(x_g^t- x_g^{t-1}))\geq 0, \\
    \hskip 0.7cm  \forall g, t=T-T_g^u+2, \cdots, T\\
	\sum_{t'=t}^{t+T_g^d-1} (1-x_g^{t'})\geq T_g^d (x_g^{t-1} - x_g^t), \\
    \hskip 0.7cm  \forall g, t=T_g^{d0}+1, \cdots, T-T_g^d+1\\
	\sum_{t'=t}^T ((1-x_g^{t'})-(x_g^{t-1} - x_g^t)) \geq 0, \\
    \hskip 0.7cm   \forall g, t = T-T_g^d+2, \cdots, T\\
	c_g^{ut} \geq C_g^u (x_g^t - x_g^{t-1}),\ \forall g,t\\
	c_g^{dt} \geq C_g^d (x_g^{t-1} - x_g^{t}),\ \forall g,t\\
	c_g^{ut}, c_g^{dt} \geq 0,\ \forall g,t\\
	x_g^t \in \{0,1\},\ \forall g,t\\
\end{cases}
\end{align}

\bibliographystyle{nonumber}

%
\end{document}